\documentclass[11pt]{article}
\usepackage{latexsym}
\usepackage{amsfonts}
\makeatletter

\def\pa#1#2{{\textstyle\partial #1\over\textstyle\partial #2}}

\def\@footnotetext#1{\insert\footins{%

\footnotesize

 \interlinepenalty\interfootnotelinepenalty

 \splittopskip\footnotesep

 \splitmaxdepth \dp\strutbox \floatingpenalty \@MM

 \hsize\columnwidth \@parboxrestore

 \edef\@currentlabel{\csname p@footnote\endcsname\@thefnmark}\@makefntext
 {\rule{\z@}{\footnotesep}\ignorespaces
#1\strut}}}

\def\abstract{\small\quotation{\hskip-\parindent\sc Abstract.}}
\def\classification{\@ifnextchar [{\@xfootnotenext}%
 {\begingroup\let\protect\noexpand
 \xdef\@thefnmark{}\endgroup
 \@footnotetext}}

\title {}

\begin{document}

\classification {{\it 1991 Mathematics Subject Classification:} Primary 14E09, 14E25; Secondary 14A10, 13B25.\\
$\dagger$) Partially supported by
CRCG Grant 25500/301/01.\\
$\ast$) Partially supported by RGC Grant Project 7126/98P.}

\begin{center}

{\bf \Large Embeddings of hypersurfaces in affine spaces}

\bigskip

{\bf Vladimir Shpilrain}$^{\dagger}$

\medskip

 and

\medskip

 {\bf Jie-Tai Yu}$^{\ast}$

\end{center}

\medskip

\begin{abstract}
\noindent    In this paper, we address the following two general problems: 
given two algebraic varieties in ${\bf C}^n$, find out whether or not 
they are (1) isomorphic; (2) equivalent under an automorphism of ${\bf C}^n$.
Although a complete solution of either of these problems is out of the
question at this time, we give here some handy and useful invariants of
isomorphic as well as of equivalent varieties. Furthermore, and more
importantly, we give a universal procedure for obtaining all possible 
algebraic varieties isomorphic to a given one, 
 and use it to construct numerous 
examples of isomorphic, but inequivalent algebraic varieties in ${\bf C}^n$.
 Among  other things, we establish the following interesting fact: 
for isomorphic hypersurfaces $\{p(x_1,...,x_n)=0\}$ and 
$\{q(x_1,...,x_n)=0\}$, the number of zeros of $grad(p)$ 
might be different from that of $grad(q)$. This implies, in particular,
that, although the fibers $\{p=0\}$ and 
$\{q=0\}$ are isomorphic, there are some other 
fibers $\{p=c\}$ and $\{q=c\}$ which are not. We construct examples
like that for any $n \ge 2$. 

\end{abstract} 

\bigskip

\noindent {\bf 1. Introduction }

\bigskip

 Let ${\bf C}[x_1,..., x_n]$ be the polynomial algebra in $n$ variables
over the field ${\bf C}$. Any collection of polynomials $p_1,...,p_m$ from
${\bf C}[x_1,..., x_n]$ determines an algebraic variety
${\{}p_i=0, i=1,...,m{\}}$ in the affine space ${\bf C}^n$.
We shall denote
this algebraic variety by $V(p_1,...,p_m)$.
 \smallskip

 We say that two algebraic varieties $V(p_1,...,p_m)$ and $V(q_1,...,q_k)$
 are {\it isomorphic} if the algebras of residue
classes ${\mathbf C}[x_1,..., x_n]/\langle p_1,...,p_m \rangle$ and
${\mathbf C}[x_1,..., x_n]/\langle q_1,...,q_k \rangle$ are
isomorphic.
Here $\langle p_1,...,p_m \rangle$ denotes the ideal of
 ${\mathbf C}[x_1,..., x_n]$ generated by $p_1,...,p_m$.
 Thus, isomorphism that we consider here is algebraic, not geometric, 
 i.e., we actually consider isomorphism of what is called  {\it affine schemes}. 
 (For example, two polynomials $x_1$ and $x_1^2$ determine the same geometric 
 set, but in our terminology, the corresponding algebraic varieties are 
not isomorphic.)
 
 On the other hand, we say that two algebraic varieties
 (or, rather,
embeddings of the same algebraic variety in ${\bf C}^n$) are
{\it equivalent} if there is an automorphism of ${\bf C}^n$ that takes one
of them onto the other. If any two embeddings of an algebraic variety in ${\bf C}^n$ are equivalent, we say that this algebraic variety
has a {\it unique embedding} in ${\bf C}^n$.

In this paper, we address two principal problems:

\smallskip

\noindent {\bf (I)} How to find out whether or not two given
algebraic varieties are isomorphic?
\smallskip

\noindent {\bf (II)} How to find out whether or not two given
algebraic varieties are equivalent? 
 \smallskip

 Both problems have received  a lot of attention over the years. For a survey on
the problem (I), we refer to \cite{Z}.  The most substantial contribution to 
the problem (II) are   \cite{AM}, \cite{AS}, \cite{Suzuki},  \cite{ZL} for equivalence 
in ${\bf C}^2$, and \cite{J}, \cite{K}, \cite{Srinivas} for equivalence 
in higher dimensions. In the latter three papers, it is shown that if 
the codimension of an algebraic variety $V$ in ${\bf C}^n$ is sufficiently
large, then $V$ has a  unique embedding in ${\bf C}^n$. Examples of non-uniquely
embedded varieties in higher dimensions are given in \cite{Harm} (for $n=5$)
  and \cite{KZ}  and \cite{Z} (for $n \ge 3$).
  
   In the present paper, we give a simple generic
procedure for constructing examples of isomorphic, but inequivalent 
varieties in any dimension. 
Furthermore, we give a very simple but efficient 
criterion for distinguishing isomorphic, but inequivalent {\it hypersurfaces}.
This  criterion allows us to show that isomorphic hypersurfaces that we  
construct, are actually inequivalent even under any {\it holomorphic} automorphism 
of the ambient space ${\bf C}^n$.

\smallskip

 However, we start by addressing the problem (I). 
Our contribution to (I) is the following Theorem 1.1. Before we give the
statement, we introduce three types of
isomorphism-preserving ``elementary"
transformations that can be applied to an arbitrary algebra of residue classes
$K[x_1,...,x_n]/R$,
where $R$ is an ideal of $K[x_1,..., x_n]$, and  $K$ an 
arbitrary ground field.
\medskip

\noindent {\bf (P1)} {\it  Introducing a new variable}: 
replace ~$K[x_1,...,x_n]/R$
~by \\ 
$K[x_1,...,x_n,y]/R+\langle y-p(x_1,...,x_n)\rangle$,
~where $p(x_1,...,x_n)$ is an arbitrary polynomial.
\smallskip

\noindent {\bf (P2)} {\it Canceling a variable} (this is the converse of (P1)):
if we have an algebra of residue classes of the form
 $K[x_1,...,x_n,y]/\langle p_1,...,p_m, q \rangle$,
where $q$ is of the form $y-p(x_1,...,x_n)$, and $p_1,...,p_m \in K[x_1,...,x_n]$,  
replace this by $K[x_1,...,x_n]/\langle p_1,...,p_m \rangle$.
\smallskip

\noindent {\bf (P3)} {\it Renaming the variables}: replace variables 
$(x_1,...,x_n)$ by $(x_{i_1},...,x_{i_n})$, where  $i_1,...,i_n$ are arbitrary distinct 
indices, not necessarily the integers  in ${\{}1,...,n{\}}$.
\smallskip

Then we have:

\medskip

\noindent {\bf Theorem 1.1.} Two algebras 
$K[x_1,..., x_n]/\langle p_1,...,p_m \rangle$ and 
$K[x_1,..., x_n]/\langle q_1,...,q_k \rangle$ are isomorphic if and only if
one can get from one of them to the other by a
sequence of transformations (P1)--(P3).

 In particular, two algebraic varieties $V(p_1,...,p_m)$ and $V(q_1,...,q_k)$
in ${\bf C}^n$ are isomorphic if and only if
one can get from the algebra 
${\mathbf C}[x_1,..., x_n]/\langle p_1,...,p_m \rangle$ to  the algebra 
${\mathbf C}[x_1,..., x_n]/\langle q_1,...,q_k \rangle$ by a
sequence of transformations  (P1)--(P3).
\medskip

This result alone does not solve the problem (I) since it does not give a hint on
how to choose the polynomials $p(x_1,...,x_n)$ in (P1).
Still, the result is quite useful because, on  one hand,
it yields   invariants (=necessary conditions of isomorphism)
 of algebraic varieties, and, on the other
hand, gives an easy practical way of constructing 
 ``exotic" examples of isomorphic algebraic varieties 
thus providing, in particular, potential examples of isomorphic, 
but inequivalent varieties.

\medskip 
   
 Theorem 1.1 also has the following interesting corollary which says 
that, although isomorphic algebraic varieties may not be 
equivalent, they are always {\it stably equivalent} (the precise 
meaning of that is clear from the statement below). Moreover, 
isomorphic varieties turn out to be  stably  equivalent 
under a {\it tame} automorphism. We call an automorphism tame
if it is a product of elementary and linear automorphisms, 
where elementary automorphisms are those that change only one
variable.

\medskip 

\noindent {\bf Corollary 1.2.} (cf. \cite{K}, \cite{Srinivas})
 If two algebraic varieties 
$~V(p_1,...,p_m)$ ~and ~$V(q_1,...,q_k)$ ~in  ${\bf C}^n$ are 
isomorphic, then the varieties ~$V(p_1,...,p_m, x_{n+1},...,x_{2n})$
~and  \\
$V(q_1,...,q_k, x_{n+1},...,x_{2n})$ in  ${\bf C}^{2n}$ are 
equivalent under a tame automorphism of ${\bf C}^{2n}$. 
\medskip

We now describe some invariants of isomorphic
varieties that can be obtained 
based on Theorem 1.1. 
Given a variety $V=V(q_1,...,q_m)$ ~(in  ${\bf C}^n$) that contains another variety 
$V'=V(r_1,...,r_s)$, consider the Jacobian matrix 
$J(q_1,...,q_m)=\left(\pa{q_i}{x_j}\right), 1\le i \le m, 1\le j \le n.$
Then, consider the image of this matrix under
the natural homomorphism from
${\mathbf C}[x_1,..., x_n]$ onto ${\mathbf C}[x_1,..., x_n]/\langle
r_1,...,r_s \rangle$.
 We denote this image
by $A_{R}(V)$, and call it the
{\it Alexander matrix} of the variety $V=V(q_1,...,q_m)$  with 
respect to the ideal $R = \langle r_1,...,r_s \rangle$.
This  resembles 
the Alexander matrix of a finitely presented
group which plays an important
role in algebraic topology -- see e.g. \cite{BZ} or  \cite{CF}. Continuing this
parallel with the classical Alexander matrix, we consider {\it elementary 
ideals} $E_k, ~k \ge 0,$ of the Jacobian  matrix $J(q_1,...,q_m)$  defined as follows: 
\smallskip

\noindent -- if  $0 < n-k \le m$, then $E_k$ is the ideal of ${\mathbf C}[x_1,..., x_n]$
 generated by all $(n-k) \times (n-k)$ minors of the matrix $J(q_1,...,q_m)$; 

\noindent -- if  $n-k > m$, then $E_k=0$;

\noindent -- if  $n-k \le 0$, then $E_k={\mathbf C}[x_1,..., x_n]$. 
\smallskip

 Denote  by $\bar E_k$ 
the image of $E_k$ in ${\mathbf C}[x_1,..., x_n]/R$. 
Then Theorem 1.1 allows us to show that $\bar E_k$ are invariants of 
isomorphic algebraic varieties. More precisely: 
\medskip

\noindent {\bf Theorem 1.3.} Let $V_1$ and $V_2$ be two  algebraic 
varieties  in ${\bf C}^n$, corresponding to the  ideals $R_1$ and $R_2$, respectively,
of the algebra ${\mathbf C}[x_1,..., x_n]$. Let $R$ be another ideal that 
contains $R_1+R_2$.  Consider 
Jacobian matrices $J_1$ and $J_2$ of the varieties $V_1$ and $V_2$.
Let $E_k^{(1)}$ and 
$E_k^{(2)}$ be elementary ideals of $J_1$ and $J_2$, respectively,
defined as above.  
Then, if the varieties $V_1$ and $V_2$ are isomorphic, one has 
$E_k^{(1)} + R=E_k^{(2)} +R,$ after possibly renaming the variables in 
$E_k^{(1)}+ R$.
\medskip

 In Section 2, we give examples of applying  Theorem 1.3 to distinguishing 
non-isomorphic algebraic varieties. Note that the condition
$E_k^{(1)} + R=E_k^{(2)} +R$ is always actually verifiable since  
  checking if two ideals of a polynomial algebra 
 are equal  can be done by using Gr\"obner 
basis algorithm (see \cite{AL}). As far as renaming the variables is
concerned, there are only finitely many ways of doing  it without changing
the whole set of variables.
\smallskip

 Then, we consider another application of Theorem 1.1, to constructing 
isomorphic, but inequivalent algebraic varieties. The most practical way of doing it
is illustrated by the following 
\smallskip

\noindent {\bf Example 1.} Let $f(x,y,z)$ be any polynomial, and  let $p(x,y,z)=x-f(x^k,y,z)$; 
$~q(x,y,z)=x-f^k(x,y,z); ~k \ge 1$. Then the varieties $V(p)$ and $V(q)$ are isomorphic.

 Indeed, the algebra $K[x,y,z]/\langle p \rangle$     
can be generated by $u=x^k, ~y$,  and   $z$ since 
in this algebra, we have $x=f(x^k,y,z)$.   These new generators $u,y,z$ are subject to the relation $u=f^k(u,y,z)$ (i.e., the minimal polynomial of $u, y$,  and   $z$ in the 
algebra $K[x,y,z]/\langle p \rangle$ is $u-f^k(u,y,z)$). Therefore, the 
varieties $V(p)$ and $V(q)$ are isomorphic.
 (In Section 3, we give a more formal and rigorous proof 
 of this fact based on our Theorem 1.1; what we have just given here, is a 
somewhat informal, ``fast" version of the same method).
\smallskip

Now that we have a procedure for constructing non-trivial examples of isomorphic 
algebraic varieties, we need invariants of {\it equivalent} varieties. It appears that in many situations, one can get away with very simple invariants. From now on, 
we are going to concentrate on {\it hypersurfaces}, i.e., on varieties of the 
form $V(p)$.
\smallskip 

 It turns out that {\it the number of zeros of the gradient} $grad(p)$ is a rather sharp 
invariant of equivalent hypersurfaces. (The fact that it {\it is} an invariant, follows 
immediately from the ``chain rule" for partial derivatives). Note that this 
number is, in general, different from 
the number of singularities of a hypersurface since a point where the gradient vanishes,
may not belong to the hypersurface. For the lack of an established name for those points, we introduce one here:

\medskip

\noindent {\bf  Definition.} Let $V(p)$ be a hypersurface in ${\bf C}^n$. A point in ${\bf C}^n$ where $grad(p)$  vanishes, is called a {\it quasi-singular} point of the hypersurface $V(p)$. 

\medskip 

Whereas the number of singular points is an intrinsic characteristic of a hypersurface (i.e.,  is invariant under an isomorphism), the number of quasi-singular points  is not, 
 and   therefore it can be used to distinguish isomorphic, but 
inequivalent hypersurfaces. 
This can be illustrated by the following  simple 

\smallskip

\noindent {\bf Example 2.} Let $p(x,y,z)=x+y+z-xyz$,  and    
$q(x,y,z)=x+y^2+yz-xyz$.  
Then    $grad(p)$ has 2 zeros, whereas $grad(q)$   vanishes nowhere. Therefore, the 
hypersurfaces $V(p)$ and   $V(q)$  are inequivalent. 
The fact that they are isomorphic, 
can be established by the same method that we used in Example 1.  
The  algebra $K[x,y,z]/\langle p \rangle$     
can be generated by $u=xy$, $y$ and   $z$ since in this algebra, $x=xyz-y-z$. 
These new generators  are subject to the relation $u=(uz-y-z)y=uyz-y^2-yz$, 
whence the isomorphism. 

\smallskip

 This example can be generalized to the following
\medskip 

\noindent {\bf Proposition 1.4.}   For any $n \ge 3$, the hypersurface 
$x_1+x_2+...+x_n-x_1x_2...x_n=0$ has at least two 
embeddings in ${\bf C}^n$ inequivalent 
  under any  algebraic or even just holomorphic automorphism of ${\bf C}^n$. 
\medskip 

 We emphasize this particular result here   because it
contrasts sharply with a result of Jelonek \cite{Jelonek} saying that
the ``$n$-cross"  $x_1x_2...x_n=0$ has a unique 
embedding in ${\bf C}^n$ for any $n \ge 3$. 
 Our Section 3 contains many other examples of isomorphic, but
inequivalent hypersurfaces.
\smallskip

 One of our  motivations for addressing these issues was the following 
well-known conjecture of Abhyankar \cite{Abhyankar} and  Sathaye \cite{S}:

\medskip 

\noindent {\bf The Embedding conjecture.}  If a hypersurface $V(p)$ in 
${\bf C}^n$ is isomorphic to the coordinate 
hyperplane  $V(x_1)$, then it is  equivalent to it. 

\medskip 

  This conjecture was proved by Abhyankar and Moh for $n=2$ \cite{AM}, and  
it remains open for  $n \ge 3$, where the prevalent opinion is that it is false, although
for  $n = 3$, there are substantial partial results in the positive 
direction \cite{Russell}. 
In order to construct a counterexample,
one needs two things: a procedure for constructing non-trivial examples 
of  hypersurfaces  isomorphic  to a  hyperplane,  and  also 
invariants of equivalent hypersurfaces, 
sharp enough to detect hypersurfaces inequivalent to  a hyperplane. 

 The procedure described above (based on Theorem 1.1) is rather flexible and 
simple. Moreover, the ``only if" part of Theorem 1.1 ensures that if a 
counterexample to the Embedding 
conjecture exists, then it can be found by our method. Hence, the main 
problem is to find 
sufficiently sharp invariants of equivalent hypersurfaces. 
The question whether or not 
the number of quasi-singular points   can possibly be useful in that 
situation, seems interesting in its own right: 

\medskip 

\noindent {\bf Problem.}  If a hypersurface $V(p)$ in ${\bf C}^n$ is 
isomorphic to  a coordinate hyperplane,  is it true that $grad(p)$ vanishes 
nowhere?

\medskip 

 If the answer to this problem is negative, that would 
also give an example of a 
polynomial some of whose fibers are 
isomorphic to $V(x_1)$, and some are not.  Indeed, suppose $X_0 \in {\bf C}^n$ 
is a   quasi-singular point of $V(p)$, and let $p(X_0)=c$. Then the 
hypersurface $V(p-c)$ has a singular point, and  therefore cannot be isomorphic 
to $V(x_1)$. 

 There is no polynomial like 
that in ${\bf C}[x_1, x_2]$ -- see \cite{AM}, but for polynomials in more 
than two variables the situation is unknown.\\

\noindent {\bf 2. Invariants of isomorphic varieties}

\bigskip

 We start with
\medskip 

\noindent {\bf  Proof of Theorem 1.1.} Since the ``if" part of the theorem 
is obvious, we proceed with the ``only if" part. 

We are going to  show  that if the algebras ~$K[X]=K[x_1,..., x_n]/\langle 
p_1,...,p_m \rangle$ ~and \\
$K[x_1,..., x_n]/\langle q_1,...,q_k \rangle$ ~are isomorphic, 
then we can take them  both to the same algebra by a sequence 
of transformations (P1)--(P3).

 We denote the former algebra by $K[X]/R$, and the latter by 
$K[X]/S$, where  $R$  and $S$ are ideals of $K[x_1,..., x_n]$
 generated by $p_1,...,p_m$ and $q_1,...,q_k$, respectively. Upon applying
 the transformation (P3), we may assume that the latter algebra is of the 
form $K[Y]/S$, where  the set $Y = {\{}y_1,..., y_n {\}}$ does not
overlap with $X$. 

 Then, by repeatedly applying the transformation (P1), we can  get 
$$K[X]/R ~\cong ~K[X \cup Y]/ 
R+ \sum_{i=1}^n \langle y_i-u_i(x_1,...,x_n)\rangle,$$ 
where $u_i$ are polynomials constructed as follows. Since $K[X]/R$ 
is isomorphic to $K[Y]/S$, there is a mapping $\varphi : 
K[X]/R \to K[Y]/S$,   such 
  that  $Ker(\varphi)=0+R$, and 
$$y_i + S = u_i(\varphi(x_1+ R),...,\varphi(x_n+ R))+ S  =
\varphi(u_i(x_1,...,x_n)+ R).$$ 

 Denote the ideal $\sum_{i=1}^n \langle y_i-u_i(x_1,...,x_n)\rangle$ by $U$. 
 Now comes the crucial point: we are going to  show  that $S \subseteq 
R + U$.

 Suppose we have an arbitrary polynomial $w(y_1,..., y_n) \in S$. 
Then  $w(y_1,..., y_n) = w((y_1-u_1)+u_1, ..., (y_n-u_n)+u_n) \equiv  
w(u_1,..., u_n) ~(mod ~U)$. On the other hand, 
$\varphi(w(u_1,..., u_n)+ R)=w(y_1,..., y_n)+ S$ ~since $y_i + S = 
\varphi(u_i+ R).$ ~Thus, $w(u_1,..., u_n)+ R \in Ker(\varphi)=0+R$, and 
therefore $w(y_1,..., y_n) \in R + U$. 
\smallskip

 Again, since $K[X]/R$ 
~is isomorphic to $K[Y]/S$, there is another mapping (the inverse of  $\varphi$)
 ~$\psi :  K[Y]/S \to K[X]/R$,   ~such 
  that  $Ker(\psi)=0+S$; 
$$x_i + R = v_i(\psi(y_1+ S),...,\psi(y_n+ S))+ R  =
\psi(v_i(y_1,...,y_n)+ S),$$ 
 and, furthermore,  $\psi(\varphi(x_i+ R))=x_i+ R$, which implies 
$$v_i(u_1, ..., u_n) \equiv  x_i ~(mod ~R).$$

 Then we have: $x_i - v_i(y_1,...,y_n) = x_i - v_i((y_1-u_1)+u_1, ..., 
 (y_n-u_n)+u_n) \equiv  x_i - v_i(u_1, ..., u_n) ~(mod ~U) \equiv 
 x_i - x_i ~(mod ~R) \equiv  0 ~(mod ~R).$
 Therefore, $x_i - v_i(y_1,...,y_n) \in R + U$.

 Thus, we have shown 
$$R + U =  R + S + U + 
\sum_{i=1}^n \langle x_i-v_i\rangle,$$
and  therefore that the algebra  $K[X]/R$ can be taken to 
$K[X \cup Y]/R + S + U + 
\sum_{i=1}^n \langle x_i-v_i\rangle$ ~by a sequence 
of transformations (P1), (P3).  By the symmetry, the algebra  
$K[Y]/S$ can be taken to the same algebra by a sequence of (P1), (P3),
 and  this completes the proof. $\Box$ 

\medskip 

\noindent {\bf Corollary 1.2} follows immediately from the proof of 
Theorem 1.1 since transformations of the type (P1) that we have used, 
are actually tame automorphisms of $K[X \cup Y]$. $\Box$

\medskip 

\noindent {\bf  Proof of Theorem 1.3} basically goes along the same 
lines as the  standard proof of invariance of Alexander  ideals  --
see e.g. \cite{CF}.  Namely, our transformations (P1) and (P2) are analogous
to Tietze transformations applied to group presentations, so that the 
invariance of $E_k+R$ can be established similarly. 

 However, there is one subtlety here. We have to show that $E_k+R$ do not
 depend on a particular choice of generators of the ideal 
 $\langle q_1,...,q_m \rangle$, i.e.,  that the ideals $E_k+R$ are well-defined. 
 Suppose ${\{}u_1,...,u_b{\}}$  is another 
  generating set of  $\langle q_1,...,q_m \rangle$, and  let 
$u_j = \sum_{i=1}^m q_i \cdot w_{ij}, ~1 \le j \le b$, where $w_{ij}$ are
some polynomials. 

 Add $b$ zero rows to  the Jacobian  matrix $J(q_1,...,q_m)$. This does not
change any of the $E_k$, which can be easily seen from the definition. 
Enumerate those zero rows somehow, and  add to the row number $j$ the 
following combination of other rows: $\sum_{i=1}^m Q_i \cdot w_{ij}$, 
where by $Q_i$ we denote the row of partial derivatives of the polynomial $q_i$.
Again, this operation does not affect any of the $E_k$. Now, modulo 
the ideal $\langle q_1,...,q_m \rangle$, the above combination of rows equals 
$\sum_{i=1}^m \partial(q_i \cdot w_{ij}) = \partial(u_j)$, 
where by $\partial(z)$   we denote the row of partial derivatives of 
$z$.  By the symmetry, we can obtain exactly
the same $(m+b) \times (m+b)$ matrix over 
${\mathbf C}[x_1,..., x_n]/\langle q_1,...,q_m \rangle$ from the Jacobian  
matrix $J(u_1,...,u_b)$, without changing any of the $E_k+R$. This completes
the proof. $\Box$ 
\medskip 

 We now illustrate Theorem 1.3 by a couple of examples. 
\medskip

\noindent {\bf Example 2.1.} Let $p=p(x,y)=x^2+xy+y^3$; $~V_1=V(p)$, 
and $q=q(x,y)=x^2+y^3$; $~V_2=V(q)$. Then $J(p)=grad(p)=(2x+y, x+3y^2)$;
$J(q)=grad(q)=(2x, 3y^2)$. Therefore, for the matrix $J(p)$ we have 
$E_1^{(1)} + R =\langle 2x+y, x+3y^2, p \rangle = 
\langle x, y \rangle $ (here we choose $R= \langle p, q \rangle$).
If we rename the variables in $E_1^{(1)}$ 
(by switching $x$ and $y$), we shall get the same ideal $\langle x, y \rangle$.

 On the other hand,  for the matrix $J(q)$ we have 
$E_1^{(2)} + R =\langle 2x, 3y^2, p, q \rangle =
\langle x, y^2 \rangle$. 
Thus, $V_1$ and   $V_2$ are not isomorphic. 
\medskip

  The point of the next example is to present a situation where our invariants
can distinguish algebraic varieties with the same topology. (We are appealing
to an intuitively clear fact that minor perturbations of coefficients of a
polynomial do not change the topology of a general fiber).  We also note that 
Makar-Limanov \cite{ML} has come up with
another way to  distinguish  varieties with the same topology 
by using purely algebraic invariants.

\medskip

\noindent {\bf Example 2.2.} Let $p=p(x,y,z)=x+yz+z^2$;
$q_k=kx^2+y^3, ~k \in {\bf C}, ~k \ne 0$.  By using our invariants, we can
show that  $~V(p,q_k)$  and  $V(p,q_m)$ are
non-isomorphic algebraic  curves in ${\bf C}^3$  if $k \ne m$.

 Let $R= \langle p, q_k, q_m \rangle$. Then the ideal $E_1^{(1)}+R$ for the
Jacobian  matrix $J(p, q_k)$ has the following (reduced) Gr\"obner basis
with respect to {\it deglex} term ordering (see \cite{AL}):
${\{}x^2, x z^2,  x + y z + z^2, 3 y^2  - 2 k x z, x y + 2 x z, z^3  + 3 x z
{\}}$. The ideal $E_1^{(2)}+R$ for the  Jacobian  matrix $J(p, q_m)$ has the
reduced Gr\"obner basis (with respect to the same term
ordering) which is obtained from the one above upon replacing 
$k$ by $m$. 

Since the reduced Gr\"obner basis with
respect to a particular term  ordering is unique, it follows that the 
ideal $E_1^{(2)}+R$ is different 
from $E_1^{(1)}+R$ if $k \ne m$ (obviously, renaming the variables cannot
change this fact), whence non-isomorphism. \\

 \noindent {\bf 3. Non-equivalent isomorphic varieties} 
\bigskip

In this section, we give various examples of inequivalent
embeddings of hypersurfaces in ${\bf C}^n$.  In all of these 
examples, we use the number of zeros of the gradient to 
establish inequivalence, as described in the Introduction. 
Therefore, in each case, we actually establish a stronger 
type of inequivalence than just being inequivalent under 
any algebraic automorphism of ${\bf C}^n$. All examples here 
are those of isomorphic hypersurfaces in ${\bf C}^n$ 
inequivalent under any {\it holomorphic} (i.e., complex 
analytic) automorphism of ${\bf C}^n$. Also, in all of these
examples, whenever we show that a given hypersurface $V(p)$ 
has inequivalent embeddings in ${\bf C}^n$, it will follow that 
$V(p)$ has inequivalent embeddings in ${\bf C}^m$ for any $m>n$ 
as well (in which case $V(p)$ is a cylindrical hypersurface).

 Finally, we note that by the examples below,  we try to 
illustrate all possible combinations of numbers of zeros of 
$grad(p)$ and $grad(q)$ 
for isomorphic hypersurfaces  $V(p)$ and $V(q)$.
That is, we have examples where $grad(p)$ vanishes nowhere, 
but $grad(q)$ has infinitely many zeros; examples where both 
$grad(p)$ and $grad(q)$ have finitely many, but different number 
of zeros, and so on. 
\medskip

 We start with the proof of Proposition 1.4, but before we 
proceed, we give a rigorous proof of isomorphism claimed in
Examples 1 and 2 in the Introduction.
\medskip

 In {\bf Example 1}, we start with the algebra 
$K[x,y,z]/\langle  p \rangle$, where $p=p(x,y,z)=
x-f(x^k,y,z)$. We are now going to give a sequence of elementary
transformations (P1)--(P3) that will bring us to the algebra 
$K[x,y,z]/\langle  q \rangle$. It will be technically more 
convenient to write those algebras of 
residue classes as ``algebras with relations", 
i.e., for example, instead of 
$K[x,y,z]/\langle p \rangle$ we shall write 
$\langle x,y,z \mid p=0\rangle$. The symbol $\cong$ below means 
``is isomorphic to". 
\smallskip 

\noindent $\langle x,y,z \mid x=f(x^k,y,z) \rangle 
~\cong ~\langle x,y,z,u \mid x=f(x^k,y,z), ~u=x^k \rangle 
= ~\langle x,y,z,u \mid x=f(u,y,z), ~u=x^k \rangle 
= ~\langle x,y,z,u \mid x=f(u,y,z), ~u=f^k(u,y,z) \rangle 
~\cong ~\langle y,z,u \mid u=f^k(u,y,z) \rangle
~\cong ~\langle x,y,z \mid x=f^k(u,y,z) \rangle$. $\Box$
\medskip

 For {\bf Example 2}, the sequence of elementary
transformations looks as follows: 
\smallskip 

\noindent $\langle x,y,z \mid x=xyz-y-z \rangle 
~\cong ~\langle x,y,z,u \mid x=xyz-y-z, ~u=xy \rangle 
= ~\langle x,y,z,u \mid x=uz-y-z, ~u=xy \rangle 
~\cong ~\langle y,z,u \mid ~u=(uz-y-z)y \rangle 
~\cong ~\langle x,y,z \mid ~x=(xz-y-z)y \rangle 
= ~\langle x,y,z \mid ~x=xyz-y^2-yz \rangle$. $\Box$
\medskip

 From now on, we are going to use a faster, slicker way 
of establishing isomorphism,  as it is done in the Introduction.
We are now ready for 
\medskip 

\noindent {\bf  Proof of  Proposition 1.4.} Let 
$p=p(x_1,...,x_n)=x_1+x_2+...+x_n-x_1x_2...x_n$. Then, since
the algebra ${\bf C}[x_1,..., x_n]/\langle  p \rangle$ can be 
generated by $x_1x_2...x_{n-1}, x_2,...,x_n$, we have,  as in 
Example 2,
the hypersurface  $V(p)$ isomorphic to $V(q)$, where
$q=x_1-(x_1x_n-x_2-...-x_n)x_2x_3...x_{n-1}=
x_1+x_2^2x_3...x_n+x_2x_3^2...x_{n-1}+...+x_2...x_n-
x_1x_2...x_n$.
\smallskip

 Now compute the gradients: $grad(p)=(1-x_2...x_n,
1-x_1x_3...x_n,..., 1-x_1x_2...x_{n-1})$. This gradient has 
$(n-1)$ zeros since, for $grad(p)=0$, one  must have 
$x_1=x_2=...=x_n$, and hence $x_1^{n-1}=1$. 

 On the other hand, $grad(q)=(1-x_2...x_n, 2x_2x_3...x_n+...
+x_3...x_{n-1}x_n-x_1x_3...x_n, ..., \\
x_2...x_{n-1}-x_1x_2...x_{n-1})$.
We are going to show that this gradient never vanishes. 

Indeed, from
the first component we see that none of the  $x_2,...,x_n$ can be
equal to 0.  Then, from the last component, we get $x_1=1$. 
Now, if we divide each component except the first and the last 
ones by $x_3...x_{n-1}$, we shall get a homogeneous system of linear 
equations in $x_2, x_3, ..., x_{n-1}$ (note that $x_n$ will cancel
out  since  $x_1=1$). The first equation is
$2x_2+x_3+...+x_{n-1}=0$,  and other equations are obtained from it
by repeatedly applying the cyclic permutation on variables. 
A system like that is easily seen to have only the trivial 
solution $x_2=x_3=...=x_{n-1}=0$, and then the first component 
of $grad(q)$ is not zero.  $\Box$

\medskip 

 One more example that seems to be worth emphasizing, 
is that of the tom Dieck -- Petrie hypersurface in ${\bf C}^n, ~n \ge 3$, 
determined 
by the polynomial $p(x,y,z)=x^2z-y^3z^2-3y^2z+2x-3y-1$ (see \cite{tom}).
 This  hypersurface has several remarkable properties -- see \cite{K2}. 
 Here we prove:
\medskip

\noindent {\bf Proposition 3.1.} Let $V(p)$  be the hypersurface in 
${\bf C}^n, ~n \ge 3$,  determined 
by  $p(x,y,z)=x^2z-y^3z^2-3y^2z+2x-3y-1$. Then there is a hypersurface 
$V(q)$ isomorphic to $V(p)$, but such  that there exists neither an algebraic 
nor a holomorphic automorphism of ${\bf C}^n$ that takes $V(p)$ onto $V(q)$.
\medskip

We note that the 
part about an algebraic automorphism was established in \cite{KZ} by 
using different, somewhat more complicated, techniques.
\medskip 

\noindent {\bf  Proof of  Proposition 3.1.} Let 
$p=p(x,y,z)=x^2z-y^3z^2-3y^2z+2x-3y-1$. Since in the algebra 
${\bf C}[x,y,z]/\langle  p \rangle$, one has $x=\frac{1}{2}y^3z^2+
\frac{3}{2}y^2z-\frac{1}{2}x^2z+\frac{3}{2}y+\frac{1}{2}$, this 
algebra can be generated by $x^2, y ,z$. Therefore,  as in 
Example 1, the hypersurface  $V(p)$ is isomorphic to $V(q)$, where
$q=x-(\frac{1}{2}y^3z^2+
\frac{3}{2}y^2z-\frac{1}{2}xz+\frac{3}{2}y+\frac{1}{2})^2$. 
\smallskip

Now compute the gradients: $grad(p)=(2+2xz, -3-3y^2z^2-6yz, 
x^2-2y^3z-3y^2)$.  This gradient has infinitely many zeros since the 
common zero locus of the three components is the same as that of
$(x-y, xz+1)$.  
\smallskip

 On the other hand, $grad(q)=(1+z\cdot Q, 
-(\frac{3}{2}y^2z^2+3yz+\frac{3}{2})Q, 
-(y^3z+\frac{3}{2}y^2-\frac{1}{2}x)Q)$, 
where $Q=Q(x,y,z)=(\frac{1}{2}y^3z^2+
\frac{3}{2}y^2z-\frac{1}{2}xz+\frac{3}{2}y+\frac{1}{2})$. 
This  gradient never vanishes. 

Indeed, the second component is $- \frac{3}{2}(yz+1)^2Q$. Since $Q\ne
0$ (otherwise, the first component would be 1), we get $yz=-1$.
 Now from the last component, we get $x=y^2$. Plug this into the 
first component and get $z=0$, which contradicts $yz=-1$. $\Box$
\medskip 

 To construct similar examples in ${\bf C}^2$ is somewhat more 
difficult since there is ``less room" there. However, we 
managed to do that as well. 
\medskip 

\noindent {\bf Example 3.2.} Let $p=p(x,y)=x^2+y^2-1$. 
We claim that the curve $V(p)$ has at least two inequivalent
embeddings in ${\bf C}^2$. (A probably difficult question is
whether or not it has a unique embedding  in ${\bf R}^2$).

 Since $x^2+y^2=(x+iy)(x-iy)$, where $i^2=-1$, the algebra 
${\bf C}[x,y]/\langle  p \rangle$ is isomorphic to 
${\bf C}[x,y]/\langle  xy-1 \rangle$. In this latter algebra,
one has $x=x^2y$, so that ${\bf C}[x,y]/\langle  xy-1 \rangle$
can be generated by $x^2, y$. Therefore, the curve  $V(p)$ 
is isomorphic to $V(q)$, where
$q=x^2y-1$. 

 Obviously, $grad(p)$ vanishes only at the origin, whereas 
$grad(q)$ has infinitely many zeros.  $\Box$
\medskip 

\noindent {\bf Example 3.3.} Let $p=p(x,y)=x-x^2-x^2y-1$.
Since in the algebra 
${\bf C}[x,y]/\langle  p \rangle$, one has $x=x^2+x^2y+1$, 
this algebra can be generated by $x^2, y$. Therefore, 
the curve  $V(p)$ is isomorphic to $V(q)$, where
$q=x-(1+x+xy)^2$. 

 Clearly, $grad(p)=(1-2x-2xy, -x^2)$ 
vanishes nowhere, whereas $grad(q)=(1-2(1+y)(1+x+xy), -2x(1+x+xy)$ 
       vanishes at the point  $(0, -\frac{1}{2})$. $\Box$
\medskip 

\noindent {\bf Example 3.4.} Let $p=p(x,y)=x-y-x^2y-x^2y^2$. 
In the algebra 
${\bf C}[x,y]/\langle  p \rangle$, one has $x=y+x^2y+x^2y^2$, 
hence this algebra can be generated by $x^2y, y$. Therefore, 
the curve  $V(p)$ is isomorphic to $V(q)$, where
$q=x-(x+y+xy)^2y$. 

 To find zeros of $grad(p)$ and $grad(q)$ is not so easy in
this example, but constructing Gr\"obner bases 
 for both  ideals  $\langle  p_x, p_y \rangle$ and
$\langle  q_x, q_y \rangle$ facilitates the process. The Gr\"obner
basis (with respect to 
the lexicographic term ordering with $y>x$) 
 for  $\langle  p_x, p_y \rangle$ turns out to be 
 ${\{}4 y^3+6^2+2y+x+2,  ~4 y^4+8y^3+4 y^2+2y+1{\}}$, 
whereas for $\langle  q_x, q_y \rangle$ it is 
${\{}36y^2+24y+8x+27, ~4y^3+4y^2+3y+1{\}}$. 

 Now we see that $grad(p)$ has 4 zeros, whereas $grad(q)$ has 
3 zeros. $\Box$ \\

\noindent {\bf Acknowledgements} 
\medskip 

 We are grateful to A. Campbell and A.  van den Essen for helpful comments,
 and to 
V.~Lin for  insightful discussions.

\noindent 
 Department of Mathematics, The City  College  of New York, New York, 
NY 10031 
\smallskip

\noindent {\it e-mail address\/}: ~shpil@groups.sci.ccny.cuny.edu \\

\noindent Department of Mathematics, The University of Hong Kong, 
Pokfulam Road, Hong Kong 

\smallskip

\noindent {\it e-mail address\/}: ~yujt@hkusua.hku.hk

\end{document}